\documentclass[12pt]{book}
\usepackage{amsmath,latexsym,amsthm}
\usepackage[dvips]{graphicx}
\usepackage{graphics}
\usepackage{amssymb}
\usepackage{color}
\usepackage{textcomp}
\usepackage{amsfonts}
\usepackage{amssymb}
\usepackage{longtable}
\usepackage{fancyhdr}
\usepackage[T2A]{fontenc}
\usepackage[utf8]{inputenc}
\usepackage[russian]{babel}
\makeatletter

\begin{document}

\begin{center}
{\Large\bf The time-fractional Airy equation on the metric graphs}
\end{center}

\begin{center}
{\sc  Rakhimov~Kamoliddin}\\
{\it National University of Uzbekistan, Tashkent, Uzbekistan\\
}
e-mail: {\tt kamoliddin\_ru@inbox.ru }

{\sc  Sobirov~Zarifboy}\\
{\it University of Geological Sciences, Uzbekistan\\
}
e-mail: {\tt sobirovzar@gmail.com }

{\sc  Jabborov~Nasridin,}\\
{\it National University of Uzbekistan\\
}
e-mail: {\tt jabborov61@mail.ru }
\end{center}

In this work we investigate Cauchy problem and initial boundary value problem for time-fractional Airy equation on the graphs with infinite and finite bonds. We studied properties of potentials for this equation and using these properties found the solutions of the considered problems. The uniqueness theorem is proved using the analogue of Gr\"{o}nwall-Bellman inequality and a-priory estimate.


%%%%%%%%%%%%%%%%%%%%%%%%%%%%%%%%%%%%%%  ДАЛЕЕ ТЕКСТ СТАТЬИ АВТОРА %%%%%%%%%%%%%%%%%%%%%%%%%%%%%%%%%%
\section{Introduction}
\label{intro}

It is known that in recent years there has been a noticeable increase of interest in the study of initial and initial-boundary value problems for equations of fractional order. This is due to the fact that the development of fractional-integral calculus gave a noticeable impetus to the study of diffusion and dispersion processes in various fields of science (see \cite{Kilbas1}, \cite{Mainardi1}, \cite{Hilfer1}, \cite{Metzler}, \cite {Phys1}).

The Schrodinger equation on the metric graph are well investigated (see \cite{Phys2}, \cite{Phys3} and references in them). Such graphs sometimes called quantum graphs. Schrodinger equation on the metric graph also investigated with Fokas unified transformation method in \cite{Xudoyberganov}.

The Airy equation on the interval were considered with Fokas unified transform method in  \cite{pelloni2004} and \cite{Himonas2019}. The potential theory for solutions of these equation were developed in \cite{Jurayev} and \cite{LC01}. The linearized Airy equation on the metric graphs was explored in the works \cite{Maqsad1}, \cite{Maqsad2}, \cite{Noja1}, \cite{Maqsad3} and \cite{Seifert1}. Also M.Cavalcante considered non linearized KdV equation in \cite{Calvante1}.

A. Pskhu was investigated properties of Airy equation with time-fractional derivative and was found its fundamental solution and investigated properties of potentials (see \cite{AVP02}). Futher, in \cite{BulletinNUUz} and \cite{Bulletin} second fundamental solution was found and was investigated the properties of the some additional potentials. Using this results they found solutions of initial and some IBVPs in infinite and finite intervals.
In this paper we consider two problems for the Airy equation with time-fractional derivative. The first one is the Cauchy problem on the open star graph with semi-infinite bonds and the second one is the initial boundary value problem (IBVP) on closed star graph with finite bonds. The solutions are found via potential's method developed in \cite{AVP02}, \cite{BulletinNUUz} and \cite{Bulletin}.
\section{Basic concepts}\label{basic}
The operator defined by the following expression is called \textit{fractional derivative (Caputo derivative)}(see \cite{AVP01})
 \begin{equation} \label{KThos}
{}_{C}D_{\eta,t}^{\alpha }g(t)=\frac{1}{\Gamma(1-\alpha)}\int\limits_{\eta}^{t}\frac{g'(\xi)}{|t-\xi|^\alpha}d\xi, \ \,\,\,\ 0<\alpha<1,
\end{equation}
where $\Gamma(x)$ is the Gamma function. And inverse of this operator is called operator of \textit{fractional integration}
\begin{equation} \label{integral}
J^\alpha_{\eta,t}g(t)=\frac{1}{\Gamma(\alpha)}\int\limits_\eta^{t}\frac{g(\xi)}{|t-\xi|^{1-\alpha}}d\xi.
\end{equation}
It is easy to show that
\begin{equation}\label{KTHxossasi}
{}_{C}D_{\eta,t}^\alpha g(t)={}_{C}D_{0 ,t-\eta}^\alpha g(t).
\end{equation}
The function \begin{equation} \label{Wright}
\phi(\lambda,\mu;z):=\sum\limits_{n=0}^\infty\frac{z^n}{n!\Gamma(\lambda n+\mu)}, \,\,\,\,\,\ \lambda>-1, \mu \in \textrm{C}
\end{equation} is called \textit{Wright function} (see \cite{FM01}).
Wright function can be represented as
$$ \phi(\lambda,\mu;z)=\frac{1}{2\pi i}\int \limits_{Ha}^{}e^{\sigma+z\sigma^{-\lambda }}\frac{d\sigma}{\sigma^\mu} $$
where the integral is taken along the Hankel contour (see \cite{AVP01}).
We have following estimate (see \cite{AVP02})
\begin{equation} \label{baho}
\left| \phi \left( -\lambda ,\mu;z \right) \right|\le C\exp \left( -\nu {{\left| z \right|}^{\frac{1}{1-\lambda }}} \right),\,\,\,C=C\left( \lambda ,\mu,\nu  \right)
\end{equation}
 where $\nu <\left( 1-\lambda  \right){{\lambda }^{\frac{\lambda}{1-\lambda }}}\cos \frac{\pi -\left| \arg z \right|}{1-\lambda},\,\,\frac{1+\lambda}{2}\pi <\left| \arg z \right|\le \pi$.
 For calculating the integral value of this function we have (see \cite{AVP02})
\begin{equation}\label{Wright_xossa}
 \int\limits_{0}^{+\infty }{\phi (-\lambda ,\mu ;az)dz}=-\frac{1}{a\Gamma (\mu +\lambda )}.
\end{equation}

\section{Cauchy problem for Airy equation}
\label{Problems1}
Let the graph has $k$ incoming and $m$ outgoing bonds. In the incoming edges, the coordinates are set from  $-\infty$ to $0$ , and in the outgoing bonds, the coordinates are set from $0$ to $+\infty$. The bonds of the graph are denoted by  $B_j, j=\overline{1,k+m}$ (Figure 1.).

%begin{figure}[hh]
%noindent\centering{\includegraphics[width=100mm]{graph.jpg}}\\
%begin{center}
%Fig. 1. Star-shaped graph
%end{figure}

On each bond of the graph, we consider the Airy equation with a fractional time derivative
\begin{equation} \label{eq}
{}_{C}D_{0,t}^{\alpha}u_j(x,t)-\frac{\partial^3}{\partial
x^3}u_j(x,t)=f_j(x,t), 0<t\leq T
\end{equation}
with boundary conditions
\begin{equation} \label{initial-C} u_j(x,0)=u_{0,j}(x), x\in \overline{B_j}, j=\overline{1,k+m}.\end{equation}

At the vertexes of the graph, we require the following gluing conditions

\begin{equation} \label{condition1}
a_ju_j(0,t)=u_1(0,t),
\end{equation}
where $j=\overline{2,k+m}.$

\begin{equation}\label{condition01}\frac{\partial u^+}{\partial
x}(0,t)=B\frac{\partial u^-}{\partial x}(0,t),
\end{equation}

\begin{equation} \label{condition2}
\sum\limits_{j=1}^{k}\frac{1}{a_j}\frac{\partial^2 u_j}{\partial
x^2}(0,t)=\sum\limits_{i=k+1}^{k+m}\frac{1}{a_i}\frac{\partial^2 u_i}{\partial
x^2}(0,t),
\end{equation}
where $a_1=1, a_j\neq0$ for $j=\overline{2,k+m},$ $u^-=(u_1,u_2,...,u_k)^T$, $u^+=(u_{k+1},u_{k+2},...,u_{k+m})^T$ and $B$- constant matrix with dimension $m\times k$.
These conditions are sometimes called Kirchhoff's conditions or the condition of conservation of flow rate at the vertex of the graph.

Functions $u_{0,j}(x), \ j=\overline{1,k+m},$ satisfy the bonding conditions (\ref{condition1}) -- (\ref{condition2}).

We will construct a regular solution of the equation (\ref{eq}) on the graph under consideration, satisfying the conditions (\ref{initial-C}) -- (\ref{condition2}), which tends to zero at $x\to \pm\infty$.

\subsection{Uniqueness of solution}\label{yagona_yechim1}
\textbf{Theorem 1.}
\emph{Let $B^{T}B-I_k$ positive defined matrix. Then problem (\ref{eq})---(\ref{condition2}) has at most one solution.}

{\sf Proof.}
In respect the following inequality \cite{Ali1}
$$\int\limits_{a}^{b}v{}_{C}D_{0,t}^{\alpha}vdx\geq\frac{1}{2}{}_{C}D_{0,t}^{\alpha}\int\limits_{a}^{b}v^2dx,$$
and relations
$$\int\limits_{a}^{b}u_{xxx}udx=uu_{xx}|_{a}^{b}-\frac{1}{2}u_{x}^{2}|_{a}^{b},$$
by using Cauchy's inequality and conditions
(\ref{initial-C}) --- (\ref{condition2}) we have

$${}_{C}D_{0,t}^{\alpha}||u||_0^2\leq(u^-)^T(I_k-B^TB)(u^-)+2||u||_0||f||_0\leq2||u||_0||f||_0\leq||u||_0^2+||f||_0^2,$$
where
$$\|u\|_0^2=\sum\limits_{j=1}^{k+m}\int\limits_{B_j}{u_j^2dx},$$
$u=(u_1,u_2,...,u_{k+m}).$

Using the analogue of Gr\"{o}nwall's inequality \cite{Ali1} we obtain the following a priori estimate from the last inequality \begin{equation}\label{aprior3}
||u||_0^2\leq||u_0||^2E_\alpha(2t^\alpha)+\Gamma(\alpha)E_{\alpha,\alpha}(2t^\alpha){}_{C}D_{0,t}^\alpha||f||_0^2.
\end{equation}

The proof of the theorem follows from (\ref{aprior3}).

\subsection{Fundamental solutions}
\label{Fundamental_yechim}
We construct the solution of the problem by the potential method. At the first, we need to get a special solution to the equation (\ref{eq}),
called fundamental solutions. A fundamental solution for the equation was found in the following form in the work \cite {AVP02}

\begin{equation} \label{fundamental1}
G_\alpha^{2\alpha/3}(x,t)=\frac{1}{3t^{1-2\alpha/3}}\left\{
    \begin{array}{ll}
    \phi(-\alpha/3,2\alpha/3;\frac{x}{t^{\alpha/3}}), &  x<0,\\
    -2\textrm{Re}[e^{{2\pi i}/3}\phi(-\alpha/3,2\alpha/3;e^{{2\pi i}/3}\frac{x}{t^{\alpha/3}})],&  x>0.
    \end{array}\right.
\end{equation}
Second fundamental solution can be written on the following form using paper \cite{Bulletin}
\begin{equation} \label{fundamental2}
V_\alpha^{2\alpha/3}(x,t)=\frac{1}{3t^{1-2\alpha/3}}\textrm{Im}[e^{{2\pi i}/3}\phi(-\alpha/3,2\alpha/3;e^{{2\pi i}/3}\frac{x}{t^{\alpha/3}})], x>0.
\end{equation}
These functions have following properties (see \cite{AVP02})
\begin{equation}\label{G-V}
{}_{C}D_{0,t}^{\nu}G_\sigma^\mu(x,t)=G_{\sigma}^{\mu-\nu}(x,t), \,\,\,\,\, \frac{\partial^3}{\partial x^3}G_\sigma^\mu(x,t)=G_\sigma^{\mu-\sigma}(x,t)
\end{equation}
and following estimate
\begin{equation}\label{G-V2}
|{}_{C}D_{0,t}^{\nu}G_\sigma^\mu(x,t)|\leq Cx^{-\theta}t^{\mu+\theta\sigma/3-1}
\end{equation}
where
$$\theta\geq\left\{ \begin{array}{cc}
& 0, \,\,\,\,\,\,\,\ (-\mu )\notin\mathbf{{N}_0} \\
& 1, \,\,\,\,\,\,\,\ (-\mu )\in \mathbf{{N}}_{0}. \\
\end{array} \right.$$
Using these functions we define following functions, called potentials
$$w_1(x,t)=\int_0^t{G_{\alpha}^{2\alpha/3}(x-a,t-\eta )\tau_1}(\eta )d\eta, \,\,\,\,\,\ w_2(x,t)=\int_0^t V_\alpha^{2\alpha/3}(x-a,t-\eta )\tau_2(\eta )d\eta,$$
$$w_3(x,t)=\int_0^t\frac{\partial^2}{\partial x^2}G_{\alpha^2{\alpha/3}}(x-a,t-\eta)\tau_3(\eta) d\eta, $$
$$w_4(x,t)=\int_0^t\frac{\partial^2}{\partial x^2}V_\alpha^{2\alpha/3}(x-a,t-\eta )\tau_4(\eta) d\eta ,$$
$$w_5(x,t)=\int_a^b G_\alpha^{2\alpha/3}(x-\xi,t )\tau_5(\xi)d\xi $$
and
$$w_6(x,t)=\int_0^t\int_a^b G_\alpha^{2\alpha/3}(x-\xi,t-\eta )f(\xi,\eta)d\xi d\eta.$$
\medskip
Let us show some properties of these functions through the following lemmas.

\textbf{Lemma 1.}
\emph{Let functions} $\tau_k(t),\ k=1,2$ \emph{are continuous and bounded on $(0;+\infty)$. Then}

1. \emph{Functions} $w_1(x,t)$ \emph{and} $w_2(x,t)$ \emph{are solutions of the equation} ${}_{C}D_{0,t}^\alpha u_j(x,t)-\frac{\partial^3 u_j(x,t)}{\partial x^3}=0$;

2. \emph{For the functions} $w_1(x,t)$ \emph{and} $w_2(x,t)$ \emph{hold relations} $$\lim_{t\to 0} w_k(x,t)=0, k=1,2.$$

\medskip

\textbf{Lemma 2.}
\emph{Let} $\tau_3(\eta), \tau_4(\eta) \in CVL(0,h).$ \emph{Then}
$$\lim_{x\to a-0} w_3(x,t)=\frac{1}{3}\tau_3(t),$$
$$\lim_{x\to a+0} w_3(x,t)=-\frac{2}{3}\tau_3(t),$$
$$\lim_{x\to a+0} w_4(x,t)=0.$$

The proofs of these lemmas are given in the work \cite{Bulletin}.
\medskip

\textbf{Lemma 3.}
\emph{Let} $\tau_5(x) \in C[a,b].$ \emph{Then the function $w_5(x,t)$ is the fundamental solution for equation (\ref{eq}) and}
$$
\lim_{t\to 0} {}_{C}D_{0,t}^{\alpha-1}w_5(x,t)=\tau_5(x).
$$

\medskip
\verb"Proof." Let show that the function $w_5(x,t)$ is the fundamental solution for equation (\ref{eq}).
Using relations (\ref{G-V}) we get
$${}_{C}D_{0,t}^{\alpha}w_5(x,t)=\int_a^b{}_{C}D_{0,t}^{\alpha}G_\alpha^{2\alpha /3}(x-\xi,t) \tau_5(\xi)d\xi=\int_a^b G_\alpha^{-\alpha /3}(x-\xi,t) \tau_5(\xi)  d\xi$$
and
$$\frac{{\partial}^{3}}{\partial x^3}w_5(x,t)=\int_a^b\frac{\partial^3}{\partial x^3} G_\alpha^{2\alpha /3}(x-\xi ,t)\tau_5(\xi)d\xi=\int_a^b G_\alpha^{-\alpha /3}( x-\xi ,t )\tau_5(\xi)d\xi.$$
Comparing these equalities we obtain that the the function $w_5(x,t)$ is the fundamental solution for equation (\ref{eq}).

Let find
$${}_{C}D_{0,t}^{\alpha -1}w_5(x,t)=\int_a^b{}_{C}D_{0,t}^{\alpha -1}G_\alpha^{2\alpha /3}( x-\xi ,t)\tau_5(\xi)d\xi=\int_a^b G_\alpha^{1-\alpha /3}(x-\xi,t)\tau_5( \xi)d\xi.$$
Using the inequality (\ref{G-V2}) we have following estimate
$$|{}_{C}D_{0,t}^{\alpha -1}w_5(x,t)|=\left| \int_a^b G_\alpha^{1-\alpha /3}(x-\xi ,t)\tau_5(\xi)d\xi \right|\le $$
$$\le\left| \max_{a\le x\le b}\tau_5(x)\int_a^b C|x-\xi|^{-\theta}t^{(1-\theta)\frac{\alpha}{3}}d\xi \right|,\,\,$$
where $1>\theta \ge 0.$
It shows that the integral form converges. Replacing $\frac{x-\xi }{{{t}^{\alpha /3}}}$ to $y$ and taking into account following computation
$$\int_{-\infty }^\infty g_\alpha(y)dy=\int_{-\infty }^\infty t^{\alpha /3}G_\alpha^{1-\alpha /3}(yt^{\alpha /3},t)dy=$$
$$=t^{\alpha /3}\int_{-\infty}^0 G_\alpha^{1-\alpha /3}(yt^{\alpha /3},t)dy+t^{\alpha /3}\int_0^{\infty}G_\alpha^{1-\alpha /3}(y t^{\alpha /3},t)dy=$$
$$ =t^{\alpha /3}\int_{-\infty}^0\frac{1}{3t^{\alpha /3}}\phi(-\frac{\alpha }{3},1-\frac{\alpha }{3};y)dy-$$
$$-2t^{\alpha /3}\textrm{Re}\left[ e^{2\pi i/3}\int_0^\infty \frac{1}{3t^{\alpha /3}}\phi(-\frac{\alpha }{3},1-\frac{\alpha }{3};e^{2\pi i/3}y)dy\right]=$$
$$ =\frac{1}{3}\int_{-\infty }^0\phi\left(-\frac{\alpha }{3},1-\frac{\alpha }{3};y \right)dy-2 \textrm{Re}\left[e^{2\pi i/3}\frac{1}{3}\int_{0}^{\infty}\phi \left( -\frac{\alpha }{3},1-\frac{\alpha }{3};e^{2\pi i/3}y \right)dy \right]=$$
$$ =\frac{1}{3}\left(\frac{1}{\Gamma(1-\alpha /3+\alpha /3)}-2 \textrm{Re}\left[ -e^{2\pi i/3}\frac{1}{e^{2\pi i/3}\Gamma(1-\alpha /3+\alpha /3)} \right] \right)=1 $$
we obtain
$$\lim_{t\to 0}{}_{C}D_{0,t}^{\alpha -1}w_5(x,t)=\lim_{t\to 0}\int_a^b G_\alpha^{1-\alpha /3}(x-\xi ,t)\tau_5(\xi )d\xi=$$
$$=\lim_{t\to 0}\int_{\frac{x-b}{t^{\alpha /3}}}^{\frac{x-a}{t^{\alpha /3}}}t^{\alpha /3}G_\alpha^{1-\alpha /3}(yt^{\alpha /3},t)\tau_5(x-t^{\alpha /3}y)dy=$$
$$=\frac{1}{3}\lim_{t\to 0}\int_{\frac{x-b}{{{t}^{\alpha /3}}}}^{\frac{x-a}{{{t}^{\alpha /3}}}}{{{g}_{\alpha }}(y)\tau_5(x-{{t}^{\alpha /3}}y)dy}=\frac{\tau_5(x)}{3}\int_{-\infty}^{+\infty}g_{\alpha}(y)dy=\tau_5(x).$$
The lemma is proved.
\medskip

\textbf{Lemma 4.}
\emph{The equation} ${}_{C}D_{0,t}^{\alpha}u(x,t)-\frac{\partial^3}{\partial x^3}u(x,t)=f(x,t)$ \emph{with initial condition}
$${}_{C}D_{0,t}^{\alpha-1}u(x,t)|_{t=0}=0$$ \emph{has a solution in the form}
$$
w_6(x,t)=\int_{0}^{t}d\eta\int_{a}^{b}G_{\alpha}^{2\alpha/3}(x-\xi,t-\eta)f(\xi,\eta)d\xi.
$$

\verb"Proof."
Using the results of \cite{AVP02} it is easy to show that the solution of Cauchy problem for the homogeneous equation
${}_{C}D_{0,t}^\alpha v(x,t)-\frac{\partial^3}{\partial x^3}v(x,t)=0$ with initial condition $v(x,0)=v_0(x)$ is
$$v(x,t)={}_{C}D_{0,t}^{\alpha-1}\int_{a}^{b}G_{\alpha}^{2\alpha/3}(x-\xi,t)v_0(\xi)d\xi.$$
Let calculate the derivatives for function $w_6(x,t)$
$${}_{C}D_{0,t}^{\alpha }w_6(x,t)=\frac{d}{dt}\int\limits_0^t d\eta \int\limits_a^b{}_{C}D_{\eta ,t}^{\alpha -1}G_\alpha^{2\alpha /3}(x-\xi ,t-\eta )f(\xi ,\eta )d\xi=$$
$$=\lim_{\eta \to t}\int\limits_a^b{}_{C}D_{\eta ,t}^{\alpha -1}G_\alpha^{2\alpha /3}(x-\xi ,t-\eta )f(\xi ,\eta )d\xi+$$
$$+\int\limits_0^t d\eta \int\limits_a^b\frac{d}{dt}{}_{C}D_{\eta ,t}^{\alpha -1}G_\alpha^{2\alpha /3}(x-\xi ,t-\eta )f(\xi ,\eta )d\xi=I_1+I_2.$$
Taking into account (\ref{KTHxossasi}) and relation (\ref{G-V}) we get
$$\int\limits_a^b{}_{C}D_{\eta ,t}^{\alpha -1}G_\alpha^{2\alpha /3}(x-\xi ,t-\eta )f(\xi ,\eta )d\xi=\int\limits_a^b{}_{C}D_{0,t-\eta }^{\alpha -1}G_\alpha ^{2\alpha /3}(x-\xi ,t-\eta )f(\xi,\eta)d\xi=$$
$$=\int\limits_a^bG_\alpha^{1-\alpha /3}(x-\xi,t-\eta )f(\xi,\eta )d\xi.$$
From relation (\ref{baho}) it follows that the integral $I_1$ converges uniformly. Substituting $\frac{x-\xi }{{{(t-\eta)}^{\alpha /3}}}$ to $y$ in this integral and taking into account the formula (\ref{Wright_xossa}) we get

$$I_1=\lim_{\eta \to t}\int\limits_a^b{}_{C}D_{\eta ,t}^{\alpha -1}G_{\alpha }^{2\alpha /3}(x-\xi ,t-\eta)f(\xi ,\eta)d\xi= \lim_{\eta \to t}\int\limits_a^bG_\alpha ^{1-\alpha /3}(x-\xi,t-\eta )f(\xi,\eta )d\xi= $$
$$\lim_{\eta \to t}\int\limits_{\frac{x-b}{( t-\eta)^{\alpha /3}}}^{\frac{x-a}{( t-\eta)^{\alpha /3}}}G_\alpha^{1-\alpha /3}((t-\eta)^{\alpha /3}y,t-\eta)f(x-( t-\eta)^{\alpha /3}y,\eta)( t-\eta)^{\alpha /3}dy= $$
$$ =\lim_{\eta \to t}\int\limits_{0}^{\frac{x-a}{{{( t-\eta )}^{\alpha /3}}}}{\frac{-2f(x-{{(t-\eta)}^{\alpha /3}}y,\eta){{(t-\eta)}^{\alpha /3}}}{3{{(t-\eta)}^{\alpha /3}}}\textrm{Re}\left( {{e}^{2\pi i/3}}\phi (-\frac{\alpha }{3},1-\frac{\alpha }{3};{{e}^{2\pi i/3}}y) \right)dy}+ $$
$$ +\lim_{\eta \to t}\int\limits_{\frac{x-b}{{{\left( t-\eta\right)}^{\alpha /3}}}}^{0}{\frac{1}{3{{\left( t-\eta \right)}^{\alpha /3}}}\phi \left( -\frac{\alpha }{3},1-\frac{\alpha }{3};y \right)f(x-{{\left( t-\eta \right)}^{\alpha /3}}y,\eta){{\left( t-\eta \right)}^{\alpha /3}}dy}= $$
$$=-\frac{2}{3}\lim_{\eta \to t}\int\limits_{0}^{\frac{x-a}{{{\left( t-\eta\right)}^{\alpha /3}}}}{\textrm{Re}\left( {{e}^{2\pi i/3}}\phi \left( -\frac{\alpha }{3},1-\frac{\alpha }{3};{{e}^{2\pi i/3}}y \right) \right)f(x-{{\left( t-\eta\right)}^{\alpha /3}}y,\eta)dy}+ $$
$$ +\frac{1}{3}\lim_{\eta \to t}\,\int\limits_{\frac{x-b}{{{\left( t-\eta  \right)}^{\alpha /3}}}}^{0}{\phi \left( -\frac{\alpha }{3},1-\frac{\alpha }{3};y \right)f(x-{{\left( t-\eta \right)}^{\alpha /3}}y,\eta)dy}=$$
$$=-\frac{2}{3}\textrm{Re}\left(\int\limits_{0}^{+\infty}{{e}^{2\pi i/3}}\phi \left( -\frac{\alpha }{3},1-\frac{\alpha }{3};{{e}^{2\pi i/3}}y \right) f(x,t)dy\right)+$$
$$+\frac{1}{3}\int\limits_{-\infty}^{0}{\phi \left( -\frac{\alpha }{3},1-\frac{\alpha }{3};y \right)f(x,t)dy}=-\frac{2}{3}\textrm{Re}\left(-{e}^{2\pi i/3} \frac{1}{e^{2\pi i/3} \Gamma(1)} \right) f(x,t) $$
$$+\frac{1}{3}\left( -\frac{1}{-\Gamma(1)}\right)f(x,t)=f(x,t).$$
Now we have $I_1=f(x,t)$. Furthermore we show that $I_2=\frac{\partial^3}{\partial x^3}u(x,t)$.
$$
{I}_{2}=\int\limits_{0}^{t}{d\eta \int\limits_{a}^{b}{\frac{d}{dt}{}_{C}D_{\eta ,t}^{\alpha -1}G_{\alpha }^{2\alpha /3}(x-\xi ,t-\eta)f(\xi ,\eta)d\xi }}=$$
$$=\int\limits_{0}^{t}{d\eta\int\limits_{a}^{b}{\frac{d}{dt}{}_{C}D_{0 ,t-\eta}^{\alpha -1}G_{\alpha }^{2\alpha /3}(x-\xi,t-\eta)f(\xi,\eta)d\xi }}=$$
\begin{equation}\label{I2}
=\int\limits_{0}^{t}{d\eta\int\limits_{a}^{b}{\frac{d}{dt}G_{\alpha }^{1-\alpha /3}(x-\xi,t-\eta)f(\xi,\eta)d\xi }}.
\end{equation}

For calculating $\frac{\partial^3}{\partial x^3}u(x,t)$ we use the relation (\ref{G-V}). So, we have
$$\frac{{{\partial }^{3}}}{\partial {{x}^{3}}}u\left(x,t \right)=\frac{{{\partial }^{3}}}{\partial {{x}^{3}}} \int\limits_{0}^{t}{d\eta \int\limits_{a}^{b}{G_{\alpha }^{2\alpha /3}(x-\xi ,t-\eta)f(\xi ,\eta)d\xi }}=$$
$$=\int\limits_{0}^{t}{d\eta \int\limits_{a}^{b}{\frac{{{\partial }^{3}}}{\partial {{x}^{3}}} G_{\alpha }^{2\alpha /3}(x-\xi ,t-\eta)f(\xi ,\eta)d\xi }}=\int\limits_{0}^{t}{d\eta \int\limits_{a}^{b}{G_{\alpha }^{2\alpha /3-\alpha}(x-\xi ,t-\eta)f(\xi ,\eta)d\xi }}=$$
\begin{equation}\label{3hosila}
\int\limits_{0}^{t}{d\eta \int\limits_{a}^{b}{\frac{{{\partial }^{3}}}{\partial {{x}^{3}}} G_{\alpha }^{-\alpha /3}(x-\xi ,t-\eta)f(\xi ,\eta)d\xi }}.
\end{equation}
Comparing (\ref{I2}) and (\ref{3hosila}) we get $I_2=\frac{\partial^3}{\partial x^3}u(x,t)$. The Lemma is proved.
\medskip

\subsection{Construction a solution of the problem. Potential method}
\label{mavjudlik1}

We seek the solutions of the problem in the form:

\begin{equation} \label{yechim}
u_j(x,t)=\int\limits_{0}^{t}G_{\alpha}^{2\alpha/3}(x-0,t-\tau)\varphi_j(\tau) d\tau+\int\limits_{0}^{t}V_{\alpha}^{2\alpha/3}(x-0,t-\tau)\psi_j(\tau) d\tau+R_j(x,t), \,\,\,\,\
\end{equation}
where $$R_j(x,t)=\int\limits_{0}^{\infty}G_{\alpha}^{2\alpha/3}(x-\xi,t-0)u_{0,j}(\xi)d\xi+\int\limits_{0}^{t}\int\limits_{-\infty}^{0}G_{\alpha}^{2\alpha/3}(x-\xi,t-0)f_j(\xi,\tau)d\xi d\tau,$$ $j=\overline{1,k+m}.$
Further, we use for abbreviate  $$ R(x,t)=(R_1(x,t),...,R_{k+m}(x,t))^T, R^-(x,t)=(R_1(x,t),...,R_k(x,t))^T, $$
$$R^+(x,t)=(R_{k+1}(x,t),...,R_{k+m}(x,t))^T, \varphi(t)=(\varphi_1(t),...,\varphi_{k+m}(t))^T,$$
$$\varphi^-(t)=(\varphi_1(t),...,\varphi_{k}(t))^T,\varphi^+(t)=(\varphi_{k+1},...,\varphi_{k+m})^T,\psi(t)=(\psi_1(t),...,\psi_{k+m}(t))^T, $$
$$\psi^+(t)=(\psi_1(t),...,\psi_{k}(t))^T,\psi^-(t)=(\psi_{k+1}(t),...,\psi_{k+m}(t))^T.$$

From the results of the work \cite{AVP02} and above potential properties it follows that these functions are the solution for equation (\ref{eq}) and they satisfy initial conditions (\ref{initial-C}).

Now we must satisfy gluing conditions on the vertexes of the graph.
By condition (\ref{condition1}), it follows that
$$a_jR_j(0,t)-R_1(0,t)=\int\limits_{0}^{t}\frac{1}{3\Gamma(2\alpha/3)(t-\tau)^{1-2\alpha/3}}\left(\varphi_1(\tau)-a_j\varphi_j(\tau)-\frac{\sqrt{3}a_j\psi_j(\tau)}{2}\right)d\tau,$$
where $j=\overline{2,k+m}.$
Further, by using definition of fractional integration (\ref{integral}) we obtain
$$J_{0,t}^{2\alpha/3}\left(\varphi_1(t)-a_j\varphi_j(t)+\frac{\sqrt3 a_j\psi_j(t)}{2}\right)=3(a_jR_j(0,t)-R_1(0,t)).$$
We may assume that
\begin{equation}\label{syst1}
\varphi_1(t)-a_j\varphi_j(t)+\frac{\sqrt3 a_j}{2}\psi_j(t)=3{}_{C}D_{0,t}^{2\alpha/3}(a_jR_j(0,t)-R_1(0,t)).
\end{equation}
Now, we calculate $u_x^{-}(0,t)$ and $u_x^{+}(0,t)$:
$$u_x^{-}(x,t)=R^{-}_x(x,t)+\int\limits_{0}^{t}\frac{1}{3(t-\tau)^{1-\alpha/3}}\phi(-\alpha/3,\alpha/3;\frac{x}{(t-\tau)^{\alpha/3}})\varphi^-(\tau)d\tau.$$
Further
$$u_x^{-}(0,t)=R^{-}_x(0,t)+ \int\limits_{0}^{t}\frac{1}{3(t-\tau)^{1-\alpha/3}}\cdot\frac{1}{\Gamma(\alpha/3)}\varphi^{-}(\tau)d\tau+$$
$$+\int\limits_{0}^{t}\frac{-2}{3(t-\tau)^{1-\alpha/3}}\mathrm{Im}[e^{\frac{4\pi i}{3}}\frac{1}{\Gamma(\alpha/3)}]\psi^{-}(\tau)d\tau=$$
$$=R^{-}_x(0,t)+\int\limits_{0}^{t}\frac{1}{3(t-\tau)^{1-\alpha/3}}\cdot\frac{1}{\Gamma(\alpha/3)}\varphi^{-}(\tau)d\tau.$$
It can be shown in the same way that
$$u_x^{+}(0,t)=R^{+}_x(0,t)+ \int\limits_{0}^{t}\frac{1}{3(t-\tau)^{1-\alpha/3}}\frac{1}{\Gamma(\alpha/3)}\varphi^{+}(\tau)d\tau+$$
$$+\int\limits_{0}^{t}\frac{\mathrm{Im}\left(e^{\frac{4\pi i}{3}}\frac{1}{\Gamma(\alpha/3)}\right)}{3(t-\tau)^{1-\alpha/3}}\psi^{+}(\tau)d\tau=R^{+}_x(0,t)+$$
$$+\int\limits_{0}^{t}\frac{1}{3(t-\tau)^{1-\alpha/3}}\cdot\frac{1}{\Gamma(\alpha/3)}\varphi^{+}(\tau)d\tau+\int\limits_{0}^{t}\frac{1}{3(t-\tau)^{1-\alpha/3}}\frac{{\sqrt3}}{2\Gamma(\alpha/3)}\psi^{+}(\tau)d\tau.$$
Using condition (\ref{condition01}) we obtain following relation
\begin{equation}\label{syst2}
B\varphi^{-}(t)-\varphi^{+}(t)+\frac{\sqrt3}{2}\psi^{+}(t)=3D_{0t}^{\alpha/3}(R^{+}_x(0,t)-BR^{-}_x(0,t)).
\end{equation}
From condition (\ref{condition2}), using \textbf{Lemma 2} for calculating $\frac{\partial^2}{\partial^2x}u_j(0,t)$, we get
$$\sum\limits_{j=1}^{k}\frac{1}{a_j}\frac{\partial^2}{\partial x^2} R_j(0,t)+\frac{1}{3}\sum\limits_{j=1}^{k}\frac{1}{a_j}\varphi_j(t)=\sum\limits_{j=k+1}^{k+m}\frac{1}{a_j}\frac{\partial^2}{\partial x^2} R_j(0,t)-\frac{2}{3}\sum\limits_{j=k+1}^{k+m}\frac{1}{a_j}\varphi_j(t).$$
So, we have
\begin{equation} \label{syst3}
\sum\limits_{j=1}^{k}\frac{1}{a_j}\varphi_j(t)+2\sum\limits_{j=k+1}^{k+m}\frac{1}{a_j}\varphi_j(t)=3\left(\sum\limits_{j=k+1}^{k+m}\frac{1}{a_j}\frac{\partial^2}{\partial x^2} R_j(0,t)-\sum\limits_{j=1}^{k}\frac{1}{a_j}\frac{\partial^2}{\partial x^2}R_j(0,t)\right).
\end{equation}
We obtained the system of integral equations (\ref{syst1}),(\ref{syst2}) and (\ref{syst3}) with respect to unknown functions $\Phi=(\varphi^-,\varphi^+,\psi^+)^T$
\begin{equation} \label{system}
M\Phi=h,
\end{equation}
where
\begin{equation} \label{matritsa}
M=\left( \begin{array}{cccccccccccc}
   1_{(k-1)\times 1} & -diag(a_2,...,a_k) & 0  \\
   1_{(m)\times 1} & -diag(a_{k+1},...,a_{k+m}) & \frac{\sqrt 3}{2}diag(a_{k+1},...,a_{k+m})  \\
   B & -I_m & \frac{\sqrt3}{2}I_m  \\
   \frac{1}{a_1}...\frac{1}{a_k} & \frac{2}{a_{k+1}}...\frac{2}{a_{k+m}} & 0  \\
\end{array} \right)
\end{equation} and
$$h=3\left( \begin{array}{cccccc}
  & {}_{C}D_{0,t}^{2\alpha /3}\left( {{a}_{2}}{{R}_{2}}\left( 0,t \right)-{{R}_{1}}\left( 0,t \right) \right) \\
 & {}_{C}D_{0,t}^{2\alpha /3}\left( {{a}_{3}}{{R}_{3}}\left( 0,t \right)-{{R}_{1}}\left( 0,t \right) \right) \\
 & ... \\
 & {}_{C}D_{0,t}^{2\alpha /3}\left( {{a}_{k+m}}{{R}_{k+m}}\left( 0,t \right)-{{R}_{1}}\left( 0,t \right) \right) \\
 & {}_{C}D_{0,t}^{2\alpha /3}\left( {{R}^{+}}\left( 0,t \right)-B{{R}^{-}}\left( 0,t \right) \right) \\
 & \sum\limits_{j=k+1}^{k+m}{\frac{1}{{{a}_{j}}}\frac{{{\partial }^{2}}}{\partial {{x}^{2}}}{{R}_{j}}\left( 0,t \right)}-\sum\limits_{j=1}^{k}{\frac{1}{{{a}_{j}}}\frac{{{\partial }^{2}}}{\partial {{x}^{2}}}{{R}_{j}}\left( 0,t \right)} \\
\end{array} \right)
.$$

Let us show that the main determinate of this system is not zero. Let $detM =0$. Then, the system $M\Phi=0$ has non trivial solutions
$\Phi = (\alpha_1,...,\alpha_{k+2m})$, which are not depended on $t$. (Note that $M$ is not depended on $t$). Then $$\varphi(t)=(\alpha_1, ... ,\alpha_{k+m}) \equiv\varphi_0,$$
$$\psi(t)=(0, ... ,0, \alpha_{k+m+1}, ... ,\alpha_{k+2m})\equiv\psi_0.$$

Substituting the obtained values of unknown functions into the solution form (\ref{yechim}),
we obtain the following solution to the problem in the case $f_j(x, t)=0, u_j(x, 0)=0$, $j =\overline{1, k + m}$:

$$u(x,t)=\varphi_0\int\limits_{0}^{t}G_{\alpha}^{2\alpha/3}(x,\tau) d\tau+\psi_0\int\limits_{0}^{t}V_{\alpha}^{2\alpha/3}(x,\tau) d\tau.$$

Using uniqueness theorem we get

$$\varphi_0G_{\alpha}^{2\alpha/3}(x,t) +\psi_0 V_{\alpha}^{2\alpha/3}(x,t)=0.$$

The last equality contradicts the linear independence of functions $G_{\alpha}^{2\alpha/3}(x,t)$ and $V_{\alpha}^{2\alpha/3}(x,t)$. So, $detM \not =0$ and
the system (\ref{system}) has unique solution.

So, we got following theorem.

\textbf{Theorem 2.}
\emph{Let} $B^{T}B-I_k$ \emph{positive defined matrix, the functions} $u_{0,j}(x)\in C(\bar B_j)$,
$f_j(x,t)\in C^{0,1}(\bar B_j\times[0,T])$, $j=\overline{1,k+m},$ \emph{and these function are zero while} $x\to\pm\infty$.
\emph{Then the problem (\ref{eq}) --- (\ref{condition2}) has unique solution on} $0\leq t\leq T$ \emph{on the form}
$$u(x,t)= R(x,t)+\int\limits_0^t U(x-\xi, t-\tau) M^{-1}h(\tau)d\tau,$$
\emph{where}
$$U(x,t)=\left(
\begin{array}{ccccc}
G_{\alpha}^{2\alpha/3}I_k &|& 0_{k\times m} &|&0_{k\times m}\\
0_{k\times m}&|& G_{\alpha}^{2\alpha/3}I_m &|& V_{\alpha}^{2\alpha/3}I_m
\end{array}\right).$$

\section{Boundary value problem for KdV equation}
    Now we consider the graph with $k$ incoming and $m$ outgoing bonds. In the incoming bonds coordinate is set from $L_j$ ($L_j<0$, $j=\overline{1,k}$)
to $0$ , and on the outgoing bonds the coordinates are set from $0$ to $L_i$ ($L_i>0, i=\overline{k+1,k+m}$). The bonds of the graph are denoted by
$b_j, j = \overline{1,k+m}$ (Figure 2).
%begin{figure}[hh]
%noindent\centering{\includegraphics[width=100mm]{graph2.jpg}}\\
%begin{center}
%Fig. 2. Star-shaped graph
%end{center}
%end{figure}

    Let us explore the equation (\ref{eq}) on the each bond $b_j$ ($j=\overline{1,k+m}$) of the above graph.
Let $t\in0\leq t\leq T,$ and $x\in\overline{b_j}, j=\overline{1,k+m}.$ We need to impose following initial conditions
\begin{equation}\label{condit01}
u(x,0)=u_0(x),
\end{equation}
vertex conditions
\begin{equation}\label{condit02}
Au(0,t)=0,
\end{equation}
\begin{equation}\label{condit03}
\frac{\partial}{\partial x}u^+(0,t)=B\frac{\partial}{\partial x}u^-(0,t),
\end{equation}
where $u^-=(u_1,u_2,...,u_k)^T,u^+=(u_{k+1},u_{k+2},...,u_{k+m})^T$,

$u=\left( \begin{array}{cc}
    u^+\\
    u^- \\
\end{array} \right)
$

$$A=\left( \begin{array}{cccccccccccccccccccccccccc}
   1 & -a_2 & 0 & ... & 0  \\
   1 & 0 & a_3 & ... & 0  \\
   ... & ... & ... & ... & ...  \\
   1 & 0 & 0 & ... & 0  \\
   1 & 0 & 0 & ... & -a_{k+m}  \\
\end{array} \right)$$
and $B$ is the constant matrix of dimension $m\times k$.

And one more vertex condition
\begin{equation}\label{condit04}
C^-{\frac{\partial^2u^-(x,t)}{\partial x^2}}|_{x=0}=C^+{\frac{\partial^2u^+(x,t)}{\partial x^2}}|_{x=0},
\end{equation}
where $C^-=(\frac{1}{a_1},\frac{1}{a_2},...,\frac{1}{a_k}),C^+=(\frac{1}{a_{k+1}},...,\frac{1}{a_{k+m}})$, $a_1=1$ and $a_j\neq0$ for $j=\overline{2,k+m}$.
Also following boundary conditions
\begin{equation}\label{condit05}
	u(L,t)=\varphi(t),
\end{equation}
\begin{equation}\label{condit06}
	\frac{\partial u^{-}(x,t)}{\partial x}|_{x=L^-}=\phi(t)
\end{equation}
where $\varphi=(\varphi_1,\varphi_2,...,\varphi_{k+m})^T$ and $\phi=(\phi_1,\phi_2,...,\phi_k)^T.$
\subsection{Uniqueness of solutions}

\textbf{Theorem 3.}
\emph{Let $B^{T}B-I_k$ negative defined matrix. Then problem (\ref{eq}), (\ref{condit01})---(\ref{condit06}) has at most one solution.}

{\sf Proof.} The proof of this theorem is similar with the proof of the \textbf{Theorem 1}.

\subsection{Existence of solutions}
Let
$$F^-=(F_1,...,F_{k})^T, F^+=(F_{k+1},...,F_{k+m})^T,\alpha^-=(\alpha_1,...,\alpha_{k})^T, \alpha^+=(\alpha_{k+1},...,\alpha_{k+m})^T,$$
$$\beta^-=(\beta_1,...,\beta_{k})^T, \beta^+=(\beta_{k+1},...,\beta_{k+m})^T, \gamma^-=(\gamma_1,...,\gamma_{k})^T, \gamma^+=(\gamma_{k+1},...,\gamma_{k+m})^T, $$
$$\rho^-=(\rho_1,\rho_2,...,\rho_{k})^T, \rho^+=\rho_{k+1},...,\rho_{k+m})^T,
 \alpha=\left( \begin{array}{cc}
  & \alpha^- \\
 & \alpha^+ \\
\end{array} \right),
\beta=\left( \begin{array}{cc}
  & \beta^- \\
 & \beta^+\\
\end{array} \right),
$$
$$\gamma=\left( \begin{array}{cc}
  & \gamma^- \\
 & \gamma^+ \\
\end{array} \right),
\rho=\left( \begin{array}{cc}
  & \rho^- \\
 & \rho^- \\
\end{array} \right),
F=\left( \begin{array}{cc}
  & F^- \\
 & F^+ \\
\end{array} \right)$$ and
$b_j=\left\{
\begin{array}{cc}
  & (L_j;0),\,\,\,\,j=\overline{1,k} \\
 & (0;L_j),\,\,\,\,j=\overline{k+1,k+m} \\
\end{array} \right. .$

We look for solutions in the form
$$u_j(x,t)=\int_0^tG_\alpha^{2\alpha /3}\left( x-L_j,t-\tau \right)\alpha_j(\tau)d\tau+\int_0^tV_\alpha^{2\alpha /3}\left(x-L_j,t-\tau \right)\beta_j(\tau )d\tau+$$
$$+\int_0^tG_\alpha^{2\alpha /3}\left(x-0,t-\tau \right)\gamma_j(\tau )d\tau+\int_0^tV_\alpha^{2\alpha /3}\left(x-0,t-\tau\right)\rho_j(\tau)d\tau+F_j(x,t),\,\,\,j=\overline{1,k+m} $$
where the functions $\alpha_j, \gamma_j$($j=\overline{1,k+m}$),  $\beta_j$($j=\overline{1,k}$), $\rho_j$($j=\overline{k+1,k+m}$) are unknown functions, $\rho_j(t)=0, (j=\overline{1,k})$; $\beta_i(t)=0, i=\overline{k+1,k+m}$ and
$$ F_j(x,t)=\int\limits_{b_j}{u_{0,j}(\xi){}_{C}D^{\alpha-1}_{0,t}G^{2\alpha/3}_\alpha (x-\xi,t-0)d\xi}+\int\limits_{0}^{t}\int\limits_{b_j}G_{\alpha}^{2\alpha/3}(x-\xi,t-0)f_j(\xi,\tau)d\xi d\tau.$$

From \textbf{Lemma 4} and the results of the work \cite{AVP02} it follows that these functions are the solutions for equation (\ref{eq}) and they satisfy initial conditions (\ref{condit01}).

Satisfying condition (\ref{condit02}) we have
$$ a_j\int_0^tG_\alpha^{2\alpha /3}(-L_j,t-\tau)\alpha_j(\tau)d\tau+a_j\int_0^tV_\alpha^{2\alpha /3}(-L_j,t-\tau)\beta_j(\tau)d\tau+$$
$${{a}_{j}}\int_{0}^{t}{G_{\alpha }^{2\alpha /3}\left( 0,t-\tau  \right){{\gamma }_{j}}(\tau )d\tau } +{{a}_{j}}\int_{0}^{t}{V_{\alpha }^{2\alpha /3}\left( 0,t-\tau  \right){{\rho }_{j}}(\tau )d\tau }+{{a}_{j}}{{F}_{j}}\left( 0,t \right)= $$
$$=\int_{0}^{t}{G_{\alpha }^{2\alpha /3}\left( -{{L}_{1}},t-\tau  \right){{\alpha }_{1}}(\tau )d\tau } +\int_{0}^{t}{V_{\alpha }^{2\alpha /3}\left( -{{L}_{1}},t-\tau  \right){{\beta }_{1}}(\tau )d\tau }+$$
$$+\int_{0}^{t}{G_{\alpha }^{2\alpha /3}\left( 0,t-\tau  \right){{\gamma }_{1}}(\tau )d\tau +{{F}_{1}}\left( 0,t \right)},\,\,j=\overline{2,k+m}.$$

Furthermore
$$\int_{0}^{t}\left({G_{\alpha }^{2\alpha /3}\left( -{{L}_{1}},t-\tau  \right){{\alpha }_{1}}(\tau )}+V_{\alpha }^{2\alpha /3}\left( -{{L}_{1}},t-\tau  \right){{\beta }_{1}}(\tau )\right)d\tau+$$
$$+\int_{0}^{t}{\frac{\phi (-\frac{\alpha }{3},\frac{2\alpha }{3};0)}{3{{(t-\tau )}^{1-2\alpha /3}}}{{\gamma }_{1}}(\tau )d\tau }+{{F}_{1}}\left( 0,t \right)=$$
$$={{a}_{j}}\int_{0}^{t}{G_{\alpha }^{2\alpha /3}\left( -{{L}_{j}},t-\tau  \right){{\alpha }_{j}}(\tau )d\tau }+{{a}_{j}}\int_{0}^{t}{V_{\alpha }^{2\alpha /3}\left( -{{L}_{j}},t-\tau  \right){{\beta }_{j}}(\tau )d\tau }+$$
$$+{{a}_{j}}\int_{0}^{t}{\frac{\phi (-\frac{\alpha }{3},\frac{2\alpha }{3};0)}{3{{(t-\tau )}^{1-2\alpha /3}}}{{\gamma }_{j}}(\tau )d\tau }+\mathrm{Im}\left[ {{a}_{j}}\int_{0}^{t}{\frac{{{e}^{2\pi i/3}}\phi (-\frac{\alpha }{3},\frac{2\alpha }{3};0)}{3{{(t-\tau )}^{1-2\alpha /3}}}{{\rho }_{j}}(\tau )d\tau } \right]+{{a}_{j}}{{F}_{j}}\left( 0,t \right).$$
So, we have
$${{a}_{j}}{{F}_{j}}\left( 0,t \right)-{{F}_{1}}\left( 0,t \right)=\int_{0}^{t}{\frac{\frac{\sqrt{3}{{a}_{j}}}{2}{{\rho }_{j}}(\tau )-{{a}_{j}}{{\gamma }_{j}}(\tau )+{{\gamma }_{1}}(\tau )}{3\Gamma \left( \frac{2\alpha }{3} \right){{(t-\tau )}^{1-2\alpha /3}}}d\tau }-$$
$$-{{a}_{j}}\int_{0}^{t}{G_{\alpha }^{2\alpha /3}\left( -{{L}_{j}},t-\tau  \right){{\alpha }_{j}}(\tau )d\tau }-{{a}_{j}}\int_{0}^{t}{V_{\alpha }^{2\alpha /3}\left( -{{L}_{j}},t-\tau  \right){{\beta }_{j}}(\tau )d\tau }+$$
$$+\int_{0}^{t}{G_{\alpha }^{2\alpha /3}\left( -{{L}_{1}},t-\tau  \right){{\alpha }_{1}}(\tau )d\tau }+\int_{0}^{t}{V_{\alpha }^{2\alpha /3}\left( -{{L}_{1}},t-\tau  \right){{\beta }_{1}}(\tau )d\tau }.$$
And
$${\gamma }_{1}(\tau )-{{a}_{j}}{{\gamma }_{j}}(\tau )+\frac{\sqrt{3}{{a}_{j}}}{2}{{\rho }_{j}}(\tau )=3{}_{C}D_{0,t}^{2\alpha /3}\left( {{a}_{j}}{{F}_{j}}\left( 0,t \right)-{{F}_{1}}\left( 0,t \right) \right)+$$
$$-3{{a}_{j}}{}_{C}D_{0,t}^{2\alpha /3}\left( \int_{0}^{t}{G_{\alpha }^{2\alpha /3}\left( -{{L}_{j}},t-\tau  \right){{\alpha }_{j}}(\tau )d\tau }+\int_{0}^{t}{V_{\alpha }^{2\alpha /3}\left( -{{L}_{j}},t-\tau  \right){{\beta }_{j}}(\tau )d\tau } \right)- $$
$$-3{}_{C}D_{0,t}^{2\alpha /3}\left( \int_{0}^{t}{G_{\alpha }^{2\alpha /3}\left( -{{L}_{1}},t-\tau  \right){{\alpha }_{1}}(\tau )d\tau }+\int_{0}^{t}{V_{\alpha }^{2\alpha /3}\left( -{{L}_{1}},t-\tau  \right){{\beta }_{1}}(\tau )d\tau } \right).$$

From above relation we obtain
$${\gamma }_{1}(\tau )-{{a}_{j}}{{\gamma }_{j}}(\tau )+\frac{\sqrt{3}{{a}_{j}}}{2}{{\rho }_{j}}(\tau )=3{}_{C}D_{0,t}^{2\alpha /3}\left( {{a}_{j}}{{F}_{j}}\left( 0,t \right)-{{F}_{1}}\left( 0,t \right) \right)-$$
$$-3\left( \int_{0}^{t}{G_{\alpha }^{0}\left( -{{L}_{1}},t-\tau  \right){{\alpha }_{1}}(\tau )d\tau }+\int_{0}^{t}{V_{\alpha }^{0}\left( -{{L}_{1}},t-\tau  \right){{\beta }_{1}}(\tau )d\tau } \right)+$$
\begin{equation} \label{system1}
+3{{a}_{j}}\left( \int_{0}^{t}{G_{\alpha }^{0}\left( -{{L}_{j}},t-\tau  \right){{\alpha }_{j}}(\tau )d\tau }+\int_{0}^{t}{V_{\alpha }^{0}\left( -{{L}_{j}},t-\tau  \right){{\beta }_{j}}(\tau )d\tau } \right) , \,\, j=\overline{2,k+m}.
\end{equation}

Analogously, from condition (\ref{condit03}) we get
$$B{{\gamma }^{-}}(t)-{\gamma }^{+}(t)+\frac{\sqrt{3}}{2}\rho^+ (t)=3{}_{C}D_{0,t}^{\alpha /3}\left(F_{x}^{+}(0,t)-BF_{x}^{-}(0,t) \right)-$$
$$-3\int_{0}^{t}{\left( BG_{\alpha }^{\alpha/3}\left(L^-,t-\tau  \right){{\alpha }^{-}}(\tau )+BV_{\alpha }^{\alpha/3}\left( -{{L}^{-}},t-\tau  \right){{\beta }^{-}}(\tau )\right)d\tau}+$$

\begin{equation} \label{system2}
+3\int_{0}^{t}\left(G_{\alpha }^{\alpha/3}\left( -{{L}^{+}},t-\tau  \right){{\alpha }^{+}}(\tau )\right)d\tau.
\end{equation}

Satisfying condition (\ref{condit04}) and using above Lemmas we have
$${{C}^{-}}{{\gamma }^{-}}(t)+2{{C}^{+}}{{\gamma }^{+}}(t)=3C\lim_{x\to 0}\frac{{{\partial }^{2}}}{\partial {{x}^{2}}}\int\limits_{0}^{t}{G_{\alpha }^{2\alpha /3}(x-L,t-\tau )\alpha (\tau )d\tau }+$$
\begin{equation} \label{system3}
+3C\lim_{x\to 0}\frac{{{\partial }^{2}}}{\partial {{x}^{2}}}\int\limits_{0}^{t}{V_{\alpha }^{2\alpha /3}(x-L,t-\tau )\beta (\tau )d\tau }+3C{{F}_{xx}}(0,t),
\end{equation}
where $C=(-C^-,C^+).$

Using conditions (\ref{condit05}) we have
$$J_{0,t}^{2\alpha /3}\left( {{\alpha }_{j}}(t)+\frac{\sqrt{3}}{2}{{\beta }_{j}}(t) \right)+\int_{0}^{t}{G_{\alpha }^{2\alpha /3}\left( {{L}_{j}},t-\tau  \right){{\gamma }_{j}}(\tau )d\tau }+$$
$$+\int_{0}^{t}{V_{\alpha }^{2\alpha /3}\left( {{L}_{j}},t-\tau  \right){{\rho }_{j}}(\tau )d\tau }+{{F}_{j}}\left( {{L}_{j}},t \right)={{\varphi }_{j}}\left( t \right), \,\,\, j=\overline(1,k+m).$$

Applying the properties of fractional operators, it follows that
$${{\alpha }_{j}}(t)+\frac{\sqrt{3}}{2}{{\beta }_{j}}(t)={}_{C}D_{0,t}^{2\alpha /3}\left( {{\varphi }_{j}}\left( t \right)-{{F}_{j}}\left( {{L}_{j}},t \right) \right)-$$
$$-{}_{C}D_{0,t}^{2\alpha /3}\left( \int_{0}^{t}{G_{\alpha }^{2\alpha /3}\left( {{L}_{j}},t-\tau  \right){{\gamma }_{j}}(\tau )d\tau }+\int_{0}^{t}{V_{\alpha }^{2\alpha /3}\left( {{L}_{j}},t-\tau  \right){{\rho }_{j}}(\tau )d\tau } \right),\,\,j=\overline{1,k+m}.$$

Above equations can be written in the following form
$$\alpha(t)+\frac{\sqrt{3}}{2}\beta(t)=-\int_{0}^{t}{G_{\alpha }^{0}(L,t-\tau)\gamma(\tau)d\tau}-\int_{0}^{t}{V_{\alpha }^{0}(L,t-\tau)\rho(\tau)d\tau}+$$
\begin{equation} \label{system4}
+{}_{C}D_{0,t}^{2\alpha /3}\left(\varphi(t)-F(L,t) \right).
\end{equation}
Analogously, from condition (\ref{condit06}) we have
$$\alpha^-(t)-\frac{\sqrt{3}}{2}\beta^-(t)=\int_{0}^{t}{G_{\alpha }^{0}(L^-,t-\tau){\gamma^-}(\tau)d\tau}+\int_{0}^{t}{V_{\alpha }^{0}(L^-,t-\tau)\rho^-(\tau)d\tau}+$$
\begin{equation} \label{system6}
+{}_{C}D_{0,t}^{\alpha /3}\left(\phi(t)-F^{-}_{x}(L^-,t) \right).
\end{equation}
We obtained following system of integral equations (\ref{system1}) --- (\ref{system6}) with respect to unknowns $\Lambda(t)$
\begin{equation}\label{sistemma}
Q\Lambda (t)+\int\limits_{0}^{t}{K(t-\tau )\Lambda (\tau )d}\tau =H,
\end{equation}
where $\Lambda$ is the unknown functions, $Q$ is the system's main matrix of dimension $(3k+3m)\times (3k+3m)$, $K$ is the matrix of potentials and $H$ is the matrix with elements defined using these coefficients.
Using above system the matrixes can be written on the form
$$ H=\left( \begin{array}{ccccc}
   -3A{}_{C}D_{0,t}^{2\alpha /3}F(0,t)  \\
   3{}_{C}D_{0,t}^{\alpha /3}\left( F_{x}^{+}(0,t)-BF_{x}^{-}(0,t) \right)  \\
   3C{{F}_{xx}}(0,t)  \\
   {}_{C}D_{0,t}^{2\alpha /3}\left( \varphi (t)-F(L,t) \right)  \\
   {}_{C}D_{0,t}^{\alpha /3}\left( \phi (t)-F_{x}^{-}(L^-,t) \right)  \\
\end{array} \right),
\Lambda=\left( \begin{array}{cccc}
  & \alpha \\
 & \beta \\
  & \gamma \\
   & \rho \\
\end{array} \right),
Q=
\left( \begin{array}{cccc}
   0 & M  \\
   Q_1 & 0  \\
\end{array} \right),$$
where $M$ is the matrix on the form (\ref{matritsa}),
${{Q }_{1}}=\left( \begin{array}{ccccccccc}
   {{I}_{k}} & 0 & \frac{\sqrt{3}}{2}{{I}_{k}}  \\
   0 & {{I}_{m}} & 0  \\
   {{I}_{k}} & 0 & -\frac{\sqrt{3}}{2}{{I}_{k}}  \\
\end{array} \right),$
and $K=\left( \begin{array}{cccc}
   K_1 & 0  \\
   0 & K_2  \\
\end{array} \right)$
where
$$K_1=3\left( \begin{array}{cccccc}
   -AG_{\alpha }^{0}(-L) & -AV_{\alpha }^{0}(-L)  \\
   -BG_{\alpha }^{0}({{L}^{-}})\,|\,\,G_{\alpha }^{0}({{-L}^{+}}) & V_{\alpha }^{0}({{-L}^{-}})  \\
   C\lim_{x\to 0}\frac{{{\partial }^{2}}}{\partial {{x}^{2}}}G_{\alpha }^{2\alpha /3}(x-L) & -C^-\lim_{x\to 0}\frac{{{\partial }^{2}}}{\partial {{x}^{2}}}V_{\alpha }^{2\alpha /3}(x-L)  \\
\end{array} \right),$$
$$K_2=\left( \begin{array}{cc}
   -G_{\alpha }^{0}(L) & -V_{\alpha }^{0}(L)  \\
   -G_{\alpha }^{0}({{L}^{-}})\,|\,0 & V_{\alpha }^{0}({{L}^{+}})  \\
\end{array} \right).$$

It is obvious that $det(Q)\neq0$ and elements of the matrix function $K(t,\tau)$ are absolutely integrable functions on $(0,T)$. So, the matrix integrable equation (\ref{sistemma}) has unique solution in $\left(C[0,t]\right)^{2k+m}$.

So, we got following theorem.

\textbf{Theorem 4.}
\emph{Let} $B^{T}B-I_k$ \emph{negative defined matrix, the functions} $u_{j,0}(x)\in C(\bar b_j)$,
$f_j(x,t)\in C^{0,1}(\bar b_j\times[0,T])$, $j=\overline{1,k+m}.$
\emph{Then the problem (\ref{eq}), (\ref{condit01})---(\ref{condit06}) has unique solution on} $0 \leq t\leq T$.

%%%%%%%%%%%%%%%%%%%%%%%%%%%%%%%%%%%%%%%%%%%%%%%%%%%%%%%%%%%%%%%%%%%%%%%%%%%%%%%%%%%%%%%%%%%%%%%%%%%%%%%%%

\end{document}